\documentclass[11pt]{amsart}
\usepackage{amssymb}

\newtheorem{theorem}{Theorem}[section]
\newtheorem{proposition}[theorem]{Proposition}
\newtheorem{corollary}[theorem]{Corollary}
\newtheorem{lemma}[theorem]{Lemma}
\theoremstyle{definition}
\newtheorem{definition}[theorem]{Definition}
\newtheorem{remark}[theorem]{Remark}

\newcommand{\ulim}{\mathrm u\mbox{-}\kern-2pt\varinjlim}
\newcommand{\tlim}{\mathrm t\mbox{-}\kern-2pt\varinjlim}
\newcommand{\glim}{\mathrm g\mbox{-}\kern-2pt\varinjlim}
\newcommand{\llim}{\mathrm l\mbox{-}\kern-2pt\varinjlim}
\newcommand{\lclim}{\mathrm{lc}\mbox{-}\kern-2pt\varinjlim}
\newcommand{\pglim}{\mathrm{pg}\mbox{-}\kern-2pt\varinjlim}
\newcommand{\mprod}{\operatornamewithlimits{\mbox{\Large$\mathsf M$}}}
\newcommand{\LP}{\operatornamewithlimits{\overrightarrow{\prod}}}

\newcommand{\LR}{\mathsf{LR}}
\newcommand{\RL}{\mathsf{RL}}
\newcommand{\Rs}{\mathsf{R}}
\newcommand{\Ls}{\mathsf{L}}
\newcommand{\BB}{\mathcal B}

\newcommand{\w}{\omega}
\newcommand{\cbox}{\boxdot}
\newcommand{\U}{\mathcal U}
\newcommand{\I}{\mathcal{I}}
\newcommand{\J}{\mathcal{J}}
\newcommand{\e}{\varepsilon}
\newcommand{\id}{\mathrm{id}}
\newcommand{\IN}{\mathbb N}
\newcommand{\IR}{\mathbb R}
\newcommand{\A}{\mathcal A}

\newcommand{\IK}{\mathbb K}
\newcommand{\PM}{\mathcal{PM}}

\newcommand{\dlim}{\varinjlim}

\title[The topology of uniform direct limits]{The topological structure of direct limits in the category of uniform spaces}
\author{Taras Banakh and Du\v san Repov\v s}
\thanks{This research was supported by Slovenian Research Agency grant P1-0292-0101, J1-9643-0101 and BI-UA/07-08-001.}

\address{Department of Mathematics, Ivan Franko National University of Lviv, and\newline
Instytut Matematyki, Uniwersytet Humanistyczno-Przyrodniczy Jana Kochanowskiego w Kielcach, Poland}
\email{tbanakh@yahoo.com}

\address{Faculty of Mathematics and Physics, and
Faculty of Education,
University of Ljubljana,
P. O. Box 2964,
Ljubljana, Slovenia 1001}
\email{dusan.repovs@guest.arnes.si}

\begin{document} 

\begin{abstract} Let $(X_n)_{n\in\w}$ be a sequence of uniform spaces such that each space $X_n$ is a subspace in $X_{n+1}$. We give an explicit description of the topology and uniformity of the direct limit $\ulim X_n$ of the sequence $(X_n)$ in the category of uniform spaces. This description implies that a function $f:\ulim X_n\to Y$ to a uniform space $Y$ is continuous if for every $n\in\IN$ the restriction $f|X_n$ is continuous and regular at the subset $X_{n-1}$ in the sense that for any entourages $U\in\U_Y$ and $V\in\U_X$ there is an entourage $V\in\U_X$ such that for each point $x\in B(X_{n-1},V)$ there is a point $x'\in X_{n-1}$ with $(x,x')\in V$ and $(f(x),f(x'))\in U$. Also we shall compare topologies of direct limits in various categories. 
\end{abstract}  
\maketitle

\section{Introduction}

Direct limits play an important role in various branches of mathematics, see 
\cite{Bier}, \cite{Dierolf}, \cite{Floret}, \cite{Glo03}, \cite{Glo07}, \cite{HSTH}, \cite{Yama}. In this paper we reveal a fundamental role of direct limits in the category of uniform spaces for understanding the topological structure of direct limits in related categories, in particular, the category of (locally) convex topological spaces and the category of topological groups. We shall give a simple description of the topology of the direct limits in the category of uniform spaces and shall apply this description to recognizing the topological structure of direct limits in some other categories. In \cite{BRLF} these results will be essentially used in the topological characterization of LF-spaces. 

By definition, an {\em LF-space} is the direct limit $\lclim X_n$ of a tower
$$X_0\subset X_1\subset X_2\subset \cdots$$of Fr\'echet (= locally convex complete metric linear) spaces  in the category of locally convex spaces, see \cite{Bier}, \cite{Floret}. Thus, $\lclim X_n$ is the linear space $X=\bigcup_{n\in\w}X_n$ endowed with the strongest topology that turns $X$ into a locally convex linear topological space such that the identity maps $X_n\to X$ are continuous. 

The union $X=\bigcup_{n\in\w}X_n$ endowed with the strongest topology making the identity maps $X_n\to X$, $n\in\w$, continuous is called the {\em topological direct limit} of the tower $(X_n)$ and is denoted by $\tlim X_n$. 

It follows from the definitions that the identity map $\tlim X_n\to\lclim X_n$ is continuous. In general, this map is not a homeomorphism, which means that the topology of topological direct limit $\tlim X_n$ can be strictly larger than the topology of the locally convex direct limit $\lclim X_n$, see \cite{Ba98}, \cite{TSH}, \cite{HSTH} or \cite{Yama}.

Between the topologies of topological and locally convex direct limits there is a  spectrum of direct limit topologies in categories that are intermediate between the category of topological and locally convex spaces. 

The most important examples of such categories are the categories of linear topological spaces, of topological groups and the category of uniform spaces.
The direct limits of a tower $(X_n)$ in those categories will be denoted by $\llim X_n$, 
$\glim X_n$, and $\ulim X_n$, respectively. The direct limit $\glim X_n$ (resp. $\llim X_n$) of a tower
$(X_n)_{n\in\w}$ of topological groups (resp. linear topological spaces) is the union $X=\bigcup_{n\in\w}X_n$ endowed with the strongest topology that turns $X$ into a topological group (resp. a linear topological space) and makes the identity maps $X_n\to X$ continuous. 

The direct limit $\ulim X_n$ of a tower $(X_n)_{n\in\w}$ of uniform spaces is defined in a similar fashion as the countable union $X=\bigcup_{n\in\w}X_n$ endowed with the strongest uniformity making the identity maps $X_n\to X$ uniformly continuous.

Each topological group $G$ will be considered as a uniform space endowed with the two-sided uniformity generated by the entourages
$U^{\LR}=\{(x,y)\in G^2:x\in yU\cap Uy\}$ where $U$ runs over symmetric neighborhoods of the neutral element in $G$. Since each continuous homomorphism $\varphi:G\to H$ between topological groups is uniformly continuous, we conclude that for any tower $(X_n)_{n\in\w}$ of locally convex spaces the identity maps 
$$\tlim X_n\to \ulim X_n\to \glim X_n\to\llim X_n\to\lclim X_n$$are continuous.
As we have already said the identity map $\tlim X_n\to\lclim X_n$ need not be a homeomorphism. 

Our crucial observation proved in Proposition~\ref{linear} below is that the identity map $\ulim X_n\to\lclim X_n$ is a homeomorphism and hence the topologies of the direct limits $\ulim X_n$, $\glim X_n$, $\llim X_n$ and $\lclim X_n$ on $X=\bigcup_{n\in\w} X_n$ coincide.

This allows us to reduce the study of the topological structure of LF-spaces to studying the topological structure of uniform direct limits $\ulim X_n$ of towers of uniform spaces. This approach will result in topological characterization of  LF-spaces given in \cite{BRLF}. This paper can be considered as the first step in realization of this program.

We start with an explicit description of the topology of the uniform direct limit $\ulim X_n$ of a tower $(X_n)_{n\in\w}$ of uniform spaces.

By a {\em tower} of uniform spaces we shall understand any increasing sequence
$$X_0\subset X_1\subset X_2\subset\cdots
$$of uniform spaces. 
By $\U_X$ we shall denote the uniformity of a uniform space $X$. For a point $x\in X$, and subsets $A\subset X$, $U\subset X^2$ let $B(x;U)=\{y\in X:(y,x)\in U\}$ and $B(A;U)=\bigcup_{a\in A}B(a;U)$ be the $U$-balls around $a$ and $A$, respectively.

The family of all subsets of $X^2$ has an interesting algebraic structure related to the operation $$U+V=\{(x,z)\in X^2:\mbox{$\exists y\in X$ such that $(x,y)\in U$ and $(y,z)\in V$}\}$$
for $U,V\subset X^2$. This operation is associative but not commutative. 

The so-defined addition operation allows us to multiply subsets $U\subset X^2$ by positive integers using the inductive formula: $1\cdot U=U$ and $(n+1)U=nU+U$ for $n>1$.

For a sequence $(U_i)_{n\in\w}$ of subsets of $X^2$ we put $\sum_{i\le n}U_i=U_0+\dots+U_n$ and $$\textstyle{\sum_{i\in\w}U_i}=\bigcup_{n\in\w}\textstyle{\sum_{i\le n}U_i.}$$ 

For a tower $(X_n)_{n\in\w}$ of sets and a point $x\in X=\bigcup_{n\in\w}X_n$ let $$|x|=\min\{n\in\w:x\in X_n\}$$ be the {\em height} of the point $x$ in $X$. 

The following theorem yields an explicit description of the topology of the uniform direct limit of a tower of uniform spaces.

\begin{theorem}\label{topbase} For any tower $(X_n)_{n\in\w}$ of uniform spaces the family $$\mathcal B=\big\{\textstyle{B(x;\sum_{i\ge|x|}U_i):x\in X,\;(U_i)_{i\ge|x|}\in\prod_{n\ge|x|}\U_{X_n}}\big\}$$is a base of the topology of the uniform direct limit $\ulim X_n$.
\end{theorem}

This theorem implies that the operation of taking uniform direct limits is topologically multiplicative in  the following sense. 

\begin{theorem} For towers $(X_n)_{n\in\w}$, $(Y_n)_{n\in\w}$ of uniform spaces the identity map $\id:\ulim (X_n\times Y_n)\to\ulim X_n\times \ulim Y_n$ is a homeomorphism.
\end{theorem}

\begin{proof} Since the identity maps $X_n\times Y_n\to\ulim X_n\times Y_n$ are uniformly continuous, so is the identity map  $\id:\ulim (X_n\times Y_n)\to\ulim X_n\times \ulim Y_n$. To show that this map is a homeomorphism, fix a neighborhood  $O(z)$ of a point $z=(x,y)$ in the space $\ulim (X_n\times Y_n)$. By Theorem~\ref{topbase}, we can assume that $O(z)=B\big(z;\sum_{n\ge|z|}W_n)$ for some sequence of entourages $(W_n)_{n\ge|z|}\in\prod_{n\ge|z|}\U_{X_n\times Y_n}$. For each $n>|z|$ choose entourages $U_n\in\U_{X_n}$ and $V_n\in\U_{Y_n}$ such that $$\big\{\big((x,y),(x',y')\big):(x,x')\in U_n,\;(y,y')\in V_n\big\}\subset W_n.$$ Since $B(z;W_{|z|})$ is a neighborhood of the point $z=(x,y)$ in $X_{|z|}\times Y_{|z|}$, we can find entourages $U\in\U_{X_{|z|}}$ and $V\in\U_{Y_{|z|}}$ such that $B(x;(|z|{+}1)U)\times B(y;(|z|{+}1)V)\subset B(z,U_{|z|})$. For all $n\le |z|$ put  $U_n=U\cap X_n^2\in\U_{X_n}$ and $V_n=V\cap Y_n^2\in\U_{Y_n}$ and observe that according to Theorem~\ref{topbase},
$$\textstyle{B(x;\sum_{i\ge|x|}U_i)\times B(y;\sum_{i\ge|y|}V_i)\subset B(z,\sum_{i\ge|z|}W_i)}=O(z)$$ is an open neighborhood of the point $z=(x,y)$ in the space $\ulim X_n\times \ulim Y_n$.
\end{proof}

\begin{corollary} For a uniform space $X$ and a tower $(Y_n)_{n\in\w}$ of uniform spaces the identity map $\ulim (X\times Y_n)\to X\times\ulim Y_n$ is a homeomorphism.
\end{corollary}

The multiplicativity property distinguishes uniform direct limits from topological direct limits, see \cite{HSTH}. According to \cite{Ba98}, the identity function $l_2\times\tlim \IR^n\to \tlim l_2\times\IR^n$ is discontinuous. Moreover, those spaces are not homeomorphic!

Theorem~\ref{topbase} will be deduced from the explicit description of the uniformity of $\ulim X_n$ given in Theorem~\ref{unifbase} below. 
The description is given in the terms of limits of monotone sequences of uniform pseudometrics.

A pseudometric on a uniform space $X$ is called {\em uniform} if for every $\e>0$ the set $$\{d<\e\}=\{(x,y)\in X^2:d(x,y)<\e\}$$ belongs to the uniformity of $X$. By \cite[8.1.11]{En}, the family $\PM_X$ of all uniform pseudometric on a uniform space $X$  generates the uniformity $\U_X$ of $X$ in the sense that the sets $\{d<1\}$, $d\in\PM_X$, form a base of the uniformity $\U_X$.

Let $(X_n)_{n\in\w}$ be a tower of uniform spaces. The definition of the uniform direct limit $\ulim X_n$ implies that a pseudometric $d$ on the uniform space $\ulim X_n$  is uniform if and only if for every $n\in\w$ the restriction $d|X_n^2$ is a uniform pseudometric on $X_n$.

A sequence of pseudometrics $(d_n)_{n\in\w}\in\prod_{n\in\w}\PM_{X_n}$ is defined to be {\em monotone} if $d_n\le d_{n+1}|X_n^2$ for any $n\in\w$. By the direct limit $\dlim d_n$ of a monotone sequence of pseudometrics $(d_n)_{n\in\w}$ we understand the pseudometric
on $X=\bigcup_{n\in\w}X_n$ defined by the formula
$$\dlim d_n(x,y)=\inf\Big\{\sum_{i=1}^n d_{|x_{i-1},\, x_i|}(x_{i-1},x_{i}):x=x_0,x_1,\dots,x_n=y\Big\}$$on $X$. In above formula $|x_{i-1},x_i|=\max\{|x_{i-1}|,|x_i|\}$ where $|x_i|=\min\{n\in\w:x\in X_n\}$ is  the height of point $x_i$ in $X$.

\begin{theorem}\label{unifbase} The uniformity of the uniform direct limit $\ulim X_n$ of a tower of uniform spaces $(X_n)_{n\in\w}$ is generated by the family of pseudometrics
$$\big\{\dlim d_n:(d_n)_{n\in\w}\in\prod_{n\in\w}\PM_{X_n}\mbox{ is monotone}\big\}.$$
\end{theorem}

 Theorem~\ref{topbase} will be used in the proof of a simple criterion of the continuity of maps defined on uniform direct limits. To state this criterion we need: 

\begin{definition} A function $f:X\to Y$ between two uniform spaces is called {\em regular at a subset} $A\subset X$ if for any entourages $U\in\U_Y$ and $V\in \U_X$ there is an entourage $W\in\U_X$ such that for any point $x\in B(A;W)$ there is a point $a\in A$ such that $(a,x)\in V$ and $(f(x),f(a))\in U$. 
\end{definition}

\begin{theorem}\label{continuity} A function $f:\ulim X_n\to Y$ defined on the uniform direct limit of a tower $(X_n)_{n\in\w}$ of uniform spaces and with values in a uniform space $Y$ is continuous if for every $n\in\IN$ the restriction $f|X_n$ is continuous and regular at the subset $X_{n-1}\subset X_n$.
\end{theorem}

\begin{corollary}\label{homeomorphism} A bijective function $h:\ulim X_n\to\ulim Y_n$ between the uniform direct limits of towers $(X_n)_{n\in\w}$, $(Y_n)_{n\in\w}$ of uniform spaces is a homeomorphism if for every $n\in\IN$ the restrictions $h|X_n$ and $h^{-1}|Y_n$ are continuous and regular at the subsets $X_{n-1}$ and $Y_{n-1}$, respectively.
\end{corollary}

In \cite{BRLF} this corollary will be used as a principal ingredient of a topological characterization of LF-spaces and some other spaces having the structure of uniform direct limit. Theorems~\ref{topbase}, \ref{unifbase} and \ref{continuity} will be proved in Sections~\ref{ptopbase}, \ref{punifbase}, and \ref{pcontinuity}, respectively. In Section~\ref{limits} we shall discuss the interplay between direct limits in various categories.

\section{Proof of Theorem~\ref{unifbase}}\label{punifbase} Theorem will follows from 
Lemmas~\ref{adequate} and \ref{pseudo} proved in this section.

Let $(X_n)_{n\in\w}$ be a tower of uniform spaces. 
In the Cartesian product $\prod_{n\in\w}\PM_{X_n}$ consider the subspace 
$$\mprod_{n\in\w}\PM_{X_n}=\big\{(d_n)_{n\in\w}\in\prod_{n\in\w}\PM_{X_n}: 
\mbox{$(d_n)_{n\in\w}$ is monotone}\big\}$$ and fix any monotone sequence of pseudometrics $(d_n)_{n\in\w}\in\mprod\limits_{n\in\w}\PM_{X_n}$. 

First we prove that in the definition of the limit pseudometric $d_\infty=\dlim d_n$ we can restrict ourselves to chains of points $x=x_0,x_1,\dots,x_n=y$ whose heights do not oscillate too much. We recall that  $|x,y|=\max\{|x|,|y|\}$ and  $|x|=\min\{n\in\w:x\in X_n\}$ for points $x,y\in X=\bigcup_{n\in\w}X_n$.

\begin{lemma}\label{modification} For any points $x,y\in X$ and $\e>0$ there is a chain of points $x=x_0,x_1,\dots,x_n=y$ such that \begin{equation}\label{eq1}
\sum_{i=1}^nd_{|x_{i-1},x_i|}(x_{i-1},x_{i})<\dlim d_n(x,y)+\e
\end{equation}
and $|x_i|<\max\{|x_{i-1}|,|x_{i+1}|\}$ for all $0<i<n$. The latter condition implies that $$|x_0|>|x_1|>\dots>|x_{s}|\le|x_{s+1}|<\dots<|x_n|$$ for some $0\le s<n$.
\end{lemma}

\begin{proof} Let us show that any sequence $x=x_0,x_1,\dots,x_n=y$ satisfying (\ref{eq1}) and having the smallest possible length $n$ has the desired property. Indeed, assume that $|x_i|\ge\max\{|x_{i-1}|,|x_{i+1}|\}$ for some $0<i<n$.

Since 
$$
d_{|x_{i{-}1}, x_i|}(x_{i{-}1},x_{i})+
d_{|x_i, x_{i{+}1}|}(x_{i},x_{i{+}1})=d_{|x_{i}|}(x_{i{-}1},x_{i})+
d_{|x_i|}(x_{i},x_{i{+}1})\ge
 d_{|x_i|}(x_{i{-}1},x_{i{+}1})\ge d_{|x_{i{-}1},x_{i{+}1}|}(x_{i{-}1},x_{i{+}1})
$$
deleting the point $x_i$ from the sequence $x_0,\dots,x_n$ will not enlarge the sum in (\ref{eq1}) but will diminish the length of the sequence, which contradicts the minimality of $n$.

It is easy to see that for the smallest number $s$ such that $|x_s|=\min_{i\le n}|x_i|$ we get $$|x_0|>|x_1|>\dots>|x_{s}|\le|x_{s+1}|<\dots<|x_n|.$$
\end{proof}

\begin{lemma}\label{uniform} For any monotone sequence of pseudometrics $(d_n)_{n\in\w}\in\mprod\limits_{n\in\w}\PM_{X_n}$ the limit pseudometric $\dlim d_n$ on $\ulim X_n$ is uniform.
\end{lemma}

\begin{proof} The uniformity of the pseudometric $d_\infty=\dlim d_n$ is equivalent to the uniform continuity of the identity map $\ulim X_n\to (X,d_\infty)$ into the pseudometric space $(X,d_\infty)$. By the  definition of the uniform direct limit $\ulim X_n$ the uniform continuity of the identity map $\ulim X_n\to(X,d_\infty)$ is equivalent to the uniform continuity of the identity embeddings $X_n\to (X,d_\infty)$ for all $n\in\w$. For every $n\in\w$ the uniform continuity of the identity embedding $X_n\to(X,d_\infty)$ trivially follows from the uniformity of the pseudometric $d_n$ and the inequality $d_\infty|X_n^2\le d_n(x,y)$ which holds according to the definition of the pseudometric $d_\infty$.
\end{proof}

A subfamily $\A\subset\mprod\limits_{n\in\w}\PM_{X_n}$ is defined to be {\em adequate} if for any sequence of entourages $(U_n)_{n\in\w}\in\prod_{n\in\w}\U_{X_n}$ there is a sequence of pseudometrics $(d_n)_{n\in\w}\in\A$ such that $\{d_n<1\}\subset U_n$ for all $n\in\w$.

\begin{lemma}\label{adequate} The family $\mprod\limits_{n\in\w}\PM_{X_n}$ is adequate.
\end{lemma}

\begin{proof} Given a sequence of entourages $(U_n)_{n\in\w}$ we need to construct a monotone sequence of uniform pseudometrics $(d_n)_{n\in\w}\in\prod_{n\in\w}\PM_{X_n}$ such that $\{d_n<1\}\subset U_n$ for every $n\in\w$. By \cite[8.1.11]{En}, for every $k\in\w$ there is a bounded uniform pseudometric $\rho_k$ on $X_k$ such that $\{\rho_k<1\}\subset U_n$. By Isbell's Extension Theorem \cite{Is} (see also \cite[8.5.6]{En}), for every $n\ge k$ the pseudometric $\rho_k$ can be extended to a uniform pseudometric $\tilde \rho_{k,n}$ on the uniform space $X_n\supset X_k$. 

For every $n\in\w$ consider the uniform pseudometric $d_n=\sum_{k\le n}\tilde\rho_{k,n}$ and observe that $(d_n)_{n\in\w}$ is a required monotone sequence of uniform pseudometrics.
\end{proof}

\begin{lemma}\label{pseudo} For any adequate family $\A\subset\mprod\limits_{n\in\w}\PM_{X_n}$ the uniformity of $\ulim X_n$ is generated by the family of limits pseudometrics $\{\dlim d_n:(d_n)_{n\in\w}\A\}$. 
\end{lemma}

\begin{proof} By Lemma~\ref{uniform}, for each sequence $(d_n)_{n\in\w}\in\A$ the limit  pseudometric $\dlim d_n$ on $\ulim X_n$ is
uniform. 

Given an entourage $U\in\U_X$ of the diagonal of the uniform space $X=\ulim X_n$, we  need to find a monotone sequence of pseudometrics $(d_n)_{n\in\w}\in\A$ such that 
$\{\dlim d_n<1\}\subset U$.

Choose a sequence of entourages $(U_n)_{n\in\w}\in\U_X^\w$ such that $5U_0\subset U$ and $2U_{n+1}\subset U_n$ for all $n\in\w$. Since the family $\A$ is adequate, there is a monotone sequence of pseudometrics $(d_n)_{n\in\w}\in\A$ such that 
$\{d_n<1\}\subset U_n$ for all $n\in\w$. We claim that the limit pseudometric $d_\infty =\dlim d_n$ has the required property: $\{d_\infty<1\}\subset U$.

Take any points $x,y\in X$ with $d_{\infty}(x,y)<1$ and applying Lemma~\ref{modification}, find a sequence of points $x=x_0,x_1,\dots x_n=y$ such that 
\begin{equation}\label{eq2}
\sum_{i=1}^n d_{|x_{i-1},x_{i}|}(x_{i-1},x_{i})<1
\end{equation}
 and $|x_0|>\dots>|x_s|\le|x_{s+1}|<\dots<|x_n|$.

Observe that for every $i\le s$ we get $|x_{i-1}|\ge|x_i|$ and then
$d_{|x_{i-1}|}(x_{i-1},x_i)=d_{|x_{i-1},x_i|}(x_{i-1},x_i)<1$. Then choice of the 
the pseudometric $d_{|x_{i-1}|}$ guarantees that $(x_{i-1},x_i)\in \{d_{|x_{i{-}1}|}<1\}\subset U_{|x_{i-1}|}$. 

For $i>s$ we get $|x_{i-1}|\le|x_i|$ and then 
$$d_{|x_i|}(x_{i-1},x_i)=d_{|x_{i-1},x_i|}(x_{i-1},x_i)<1$$implies  $(x_{i-1},x_i)\in U_{|x_i|}$. 

Since $|x_{0}|>\dots>|x_{s}|\le |x_{s+1}|<\dots<|x_n|$, we see that $$
(x,y)\in U_{|x_0|}+\dots+U_{|x_{s-1}|}+U_{|x_{s}|}+\dots+U_{|x_{n-1}|}\subset U_{s}+\dots + U_0+U_0+U_0+U_1+\dots+U_{n-s-1}\subset 5U_0\subset U.
$$
\end{proof}

\section{Proof of Theorem~\ref{topbase}}\label{ptopbase}

Given a tower $(X_n)_{n\in\w}$ of uniform spaces we need to check that the 
 family $$\textstyle{\mathcal B=\big\{B(x;\sum_{n\ge|x|}U_n):x\in X,\;(U_n)_{n\ge|x|}\in\prod_{n\ge|x|}\U_{X_n}\big\}}$$ is a base of the topology of the uniform direct limit $\ulim X_n$.

First we prove that each set $B(x;\sum_{n\ge|x|}U_n)\in \mathcal B$ is a neighborhood of $x$ in $\ulim X_n$. 

Let $\|x\|=\min\{n\le |x|:x\in\overline{X}_n\}$ where $\overline{X}_n$ is the closure of $X_n$ in $X_{|x|}$. By definition of $\|x\|$, there is an entourage $V_{|x|}\in\U_{X_{|x|}}$ such that $B(x;2V_{|x|})\cap X_{\|x\|-1}=\emptyset$. We can take $V_{|x|}$ so small that $(|x|+2)\,V_{|x|}\subset U_{|x|}$. Put $V_k=V_{|x|}\cap X_k^2$ for all $k\le|x|$ and $V_k=U_k$ for all $k>|x|$.

By Lemma~\ref{adequate}, there is monotone sequence of uniform pseudometrics $(d_n)_{n\in\w}\in\prod_{n\in\w}\PM_{X_n}$ such that 
$\{d_n<1\}\subset V_k$ for all $k\in\w$. The limit pseudometric $d_\infty=\dlim d_n$ determines the unit ball $B_1(x)=\{y\in X:d_{\infty}(x,y)<1\}$ which is open in the space $\ulim X_n$ according to Theorem~\ref{unifbase}.

We claim that $B_1(x)\subset B(x;\sum_{n\ge|x|}U_n)$. Given any point $x'\in B_1(x)$, we need to check that $x'\in B(x;\sum_{n\ge|x|}U_n)$. Since $x\in\overline{X}_{\|x\|}$, we may find a point $x''\in B(x,V_{|x|})\cap \overline{X}_{\|x\|}$ so close to $x$ that $d(x,x'')<1-d(x,x')$. It follows from $$B(x'',V_{|x|})\cap X_{\|x\|-1}\subset B(x,2V_{|x|})\cap X_{\|x\|-1}=\emptyset$$ that $\|x''\|=|x''|=\|x\|$.

Applying Lemma~\ref{modification}, find a sequence of points $x''=x_0,x_1,\dots,x_m=x'$ such that
\begin{equation}\label{sum3}\sum_{i=1}^m d_{|x_{i-1},x_i|}(x_{i-1},x_i)<1\end{equation} and 
\begin{equation}\label{eq3}
|x_i|<|x_{i-1},x_i|=\max\{|x_{i-1}|,|x_{i+1}|\}\mbox{ for all $0<i<m$.}
\end{equation} Let $k\le m$ be the largest number such that $|x_k|=\min\{|x_i|:i\le m\}$. We claim that $k\le 1$. In the opposite case, the condition (\ref{eq3}) we would imply $|x_1|<|x_0|=\|x\|$ and then $d_{|x_0|}(x_0,x_1)<1$ by (\ref{sum3}). Since $x''=x_0,x_1\in X_{|x_0|}$, the choice of the pseudometric $d_{|x_0|}$ ensures that $(x_0,x_1)\in \{d_{|x_0|}<1\}\subset V_{|x_0|}\subset V_{|x|}$. Consequently, $x_1\in X_{|x_0|-1}\cap B(x_0;V_{|x|})\subset X_{\|x\|-1}\cap B(x'';V_{|x|})=\emptyset$, which is a contradiction.

Therefore $k\le 1$ and $\|x\|=|x_0|\le|x_1|<|x_2|<\dots<|x_m|$. It follows from (\ref{sum3}) that for every $0<i\le m$ we get $d_{|x_i|}(x_{i-1},x_i)<1$.
Since $x_{i-1},x_i\in X_{|x_i|}$, the choice of the pseudometric $d_{|x_i|}$  guarantees that $(x_{i-1},x_i)\in \{d_{|x_i|}<1\}\subset V_{|x_i|}$. Observe that the number $p=\max\{i\le m:|x_i|\le|x|\}$ does not exceed $|x|+1$. It follows from $(x,x_0)\in V_{|x|}$ and $(x_{i-1},x_i)\in V_{|x_i|}\subset V_{|x|}$, $i\le p$, that $(x,x_p)\subset (p+1)\cdot V_{|x|}\subset (|x|+2)V_{|x|}\subset U_{|x|}$. Consequently, $x_p\in B(x;U_{|x|})$. If $p=m$, then $x'=x_m\in B(x;U_{|x|})\subset B(x;\sum_{i\ge|x|}U_i)$ and we are done. If $p<m$, then can use the relation $(x_{i-1},x_i)\in V_{|x_i|}\subset U_{|x_i|}$, $p<i\le m$, in order to prove by induction that for every $i\in\{p+1,\dots,m\}$ we get
$$x_i\in B(x_p;\sum_{p<j\le i}U_{|x_j|})\subset B(x_p,\sum_{|x|<j\le|x_i|}U_j)\subset B(x;\sum_{|x|\le j\le |x_i|}U_j)\subset B(x;\sum_{j\ge|x|}U_j).$$ In particular, $x'=x_m\in B(x;\sum_{j\ge|x|}U_j)$.
This completes the proof of the inclusion $B_1(x)\subset B(x,\sum_{j\ge|x|}U_j)$. 
\smallskip

Next, we show that each set $B(x;\sum_{i\ge|x|}U_i)\in\mathcal B$ is open. Given any point $y\in B(x; \sum_{i\ge|x|}U_i)$, find the smallest number $m\in\w$ with $y\in B(x;\sum_{|x|\le i\le m}U_i)$. Since  $B(x;U_{m})$ is a neighborhood of the point $y$ in $X_m$, there is an entourage $V\in\U_{X_m}$ such that $B(y;(m+1)V)\subset B(x;U_m)$. Define a sequence of entourages $(V_i)_{i\ge|y|}\in\prod_{i\ge|y|}\U_{X_i}$ letting $V_i=V\cap X_i^2$ for $|y|\le i\le m$ and $V_i=U_i$ for $i>m$. 
It follows that $$\sum_{i\ge|y|}V_i=\sum_{|y|\le i\le m}V_i+\sum_{i>m}V_i\subset (m+1)V+\sum_{i>m}V_i$$and hence

$$
\begin{aligned}
B(y;\sum_{i\ge|y|}V_i)\subset &\,B(y;(m+1)V+\sum_{i>m}V_i)\subset B(B(y;(m{+}1)V);\sum_{i>m}V_i)\subset\\
& B(B(x;U_m);\sum_{i>m}U_i)=B(x;\sum_{i\ge m}U_i)\subset B(x;\sum_{i\ge |x|}U_i).
\end{aligned}$$ 
Since $y$ is an interior point of $B(y,\sum_{i\ge|y|}V_i)$, it is an interior point of $B(x;\sum_{i\ge|x|}U_i)$ as well.

In such a way we proved that the family $\mathcal B$ consists of open subsets of $\ulim X_n$. To show that $\mathcal B$ is a base of the topology of $X=\ulim X_n$, fix any point $x\in X$ and an entourage $U\in\U_X$. By induction find a sequence of entourages $(U_n)_{n\in\w}\in\U_X^\w$ such that $2U_0\subset U$ and $2U_{n+1}\subset U_n$ for all $n\in\w$.
For this sequence we get $\sum_{i\in\w}U_i\subset 2U_0\subset U$.
 Define a sequence of entourages $(V_i)_{i\in\w}\in\prod_{i\in\w}\U_{X_i}$ letting $V_i=U_i\cap X_i^2$ for $i\in\w$, and observe that $B(x;\sum_{i\ge|x|}V_i)\subset B(x;\sum_{i\ge|x|}U_i)\subset B(x;U)$.

\section{Proof of Theorem~\ref{continuity}}\label{pcontinuity}

Let $(X_n)_{n\in\w}$ be a tower of uniform spaces and $f:\ulim X_n\to Y$ be 
a function into a uniform space $Y$ such that for every $n\in\w$ the restriction $f|X_n$ is continuous and regular at the closed subset $X_{n-1}\subset X_n$. Let $X=\ulim X_n$.

We need to shall check the continuity of $f$ at an arbitrary point $x_0\in \ulim X_n$. 
Without loss of generality, $x_0\in X_0$. Given any entourage $U\in\U_Y$ we need to find a neighborhood $O(x_0)\subset X$ of the point $x_0$ such that $f(O(x_0))\subset B(y_0;U)$ where $y_0=f(x_0)$.

By induction construct a sequence of entourages $U_n\in\U_X$, $n\in\w$ such that  \begin{equation}\label{eq5:1}
\mbox{$3U_0\subset U$ and $2U_{n+1}\subset U_n$ for all $n\in\w$.}
\end{equation} For such a sequence we get the inclusion $\sum_{i\in\w}U_i\subset 2U_0\subset U$.

The continuity of the restriction $f|X_0$ yields an entourage $W_0\in\U_{X_0}$  such that $f(B(x_0;W_0))\subset B(y_0;U_0)$. Let us recall that for every $k\in\w$ the map $f|X_{k+1}:X_{k+1}\to Y$ is regular at the subset $X_k$. This fact can be used to construct inductively two sequences of entourages $(V_k)_{k\in\w},(W_k)_{k\in\w}\in\prod_{k\in\w}\U_{X_k}$ such that for every $k\in\w$ we get  
\begin{enumerate}
\item $X_{k-1}^2\cap 3V_{k}\subset W_{k-1}\subset V_{k-1}$, and 
\item for every $x\in B(X_{k-1},2W_k)$ there is a point $x'\in X_{k-1}$ such that $(x,x')\in V_k$ and $(f(x),f(x')\in U_k$.
\end{enumerate}

By Theorem~\ref{topbase}, the set $O(x_0)=B(x_0;\sum_{i\ge0}W_i)$ is a neighborhood of $x_0$ in $\ulim X_n$. We claim that $f(O(x_0))\subset B(f(x_0);U)$. 

Given any point $x\in B(x_0;\sum_{i\ge0}W_i)$, find $k>0$ with $x\in B(x_0;\sum_{i\le k}W_i)\subset X_k$ and put $x_k=x$. Since $x\in B(x_0;\sum_{i\le k}W_i)=B(B(x_0;\sum_{i<k}W_i);W_k)$, there is a point $x_{k-1}\in B(x_0,\sum_{i<k}W_i)$ such that $(x_k,x_{k-1})\in W_k$. Continuing by induction, we shall construct a sequence of points $x_k,x_{k-1},x_{k-2},\dots,x_1$ such that $x_i\in B(x_0;\sum_{j\le i}W_j)$ and $(x_i,x_{i-1})\in W_i$ for all $1\le i\le k$.

Let $z_k=x=x_k$. Since $z_k\in B(x_0;\sum_{i\le k}W_i)\subset B(X_{k-1};W_k)$,
there is a point $z_{k-1}\in X_{k-1}$ such that $(z_{k-1},z_k)\in V_k$ and
$(f(z_k),f(z_{k-1})\in U_k$.

It follows from $(z_{k-1},z_k)\in V_k$, $z_k=x_k$, and $(x_k,x_{k-1})\in W_k$ that 
$(z_{k-1},x_{k-1})\in V_k+W_k\subset 2V_k$ and thus  $(z_{k-1},x_{k-1})\in X_{k-1}^2\cap 2 V_k\subset W_{k-1}$. 
Then $$(z_{k-1},x_{k-2})\le 2W_{k-1}$$and hence $z_{k-1}\in B(X_{k-2};2W_{k-1})$. The choice of $W_{k-1}$ yields a point $z_{k-2}\in X_{k-2}$ such that $(z_{k-2},z_{k-1})\in V_{k-1}$ and $(f(z_{k-2}),f(z_{k-2}))\in U_{k-1}$.

For this point $z_{k-2}$ we get
$$(z_{k-2},x_{k-2})\in\{(z_{k-2},z_{k-1})\}+\{(z_{k-1},x_{k-1})\}+\{(x_{k-1},x_{k-2})\}\subset 
V_{k-1}+W_{k-1}+W_{k-1}\subset 3V_{k-1}
$$
and thus $(z_{k-2},x_{k-2})\in X_{k-2}^2\cap 3V_{k-1}\subset W_{k-2}$.
 
Continuing by induction, we will define a sequence of points $x=z_k,z_{k-1},\dots,z_0$ such that $(z_i,x_i)\in W_i$ and $(f(z_{i+1}),f(z_{i}))\in U_{i+1}$ for all $i<k$. It follows from (\ref{eq5:1}) that $(f(z_k),f(z_0))\in 2U_0.$

Since $(z_0,x_0)\in W_0$, we get $f(z_0)\in f(B(x_0;W_0))\subset B(y_0;U_0)$ and hence $$f(x)=f(z_k)\in B(f(z_0);2U_0)\subset B(B(y_0;U_0);2U_0)=B(y_0;3U_0)\subset B(y_0;U).$$

\section{Interplay between the topologies of the direct limits in various categories}\label{limits}

In this section we shall apply Theorem~\ref{topbase} to show that the topology of the direct limit in the category of uniform spaces coincides with the topology of direct limit in some other categories related to topological algebra or functional analysis. 
One of such categories in the category of abelian (more generally, SIN) topological groups and their continuous homomorphisms. 

Each topological group $G$ carries four natural uniformities compatible with the topology:
\begin{itemize}
\item[1)] the {\em left uniformity\/} $\U^{\Ls}$ generated by the entourages $U^{\Ls}=\{(x,y)\in G:x\in yU\}$,
\item[2)] the {\em right uniformity\/} $\U^{\Rs}$ generated by the entourages $U^{\Rs}=\{(x,y)\in G:x\in Uy\}$,
\item[3)] the {\em two-sided uniformity\/} $\U^{\LR}$ generated by the entourages  $U^{\LR}=\{(x,y)\in G:x\in yU\cap Uy\}$,
\item[4)] the {\em Roelcke uniformity\/} $\U^{\RL}$, generated by the entourages $U^{\RL}=\{(x,y)\in G:x\in UyU\}$,
\end{itemize}
where $U$ runs over open symmetric neighborhoods of the neutral element $e$ in the topological group $G$.

These four uniformities on $G$ coincide if and only if $G$ is a SIN-group. The latter means that $G$ has a neighborhood base at $e$, consisting of open symmetric neighborhoods $U\subset G$ that are {\em invariant} in the sense that $xUx^{-1}=U$ for all $x\in G$, see \cite{RD}.

In the sequel saying about uniform properties of topological groups we shall refer to the two-sided uniformity. In case of a SIN-group, the two-sided uniformity coincides with the other three uniformities. 

\begin{proposition}\label{group}
The uniform direct limit $\ulim X_n$ of a tower $(X_n)_{n\in\w}$ of SIN-groups is a topological group. Consequently, the identity map $\ulim X_n\to \glim X_n$ is a homeomorphism.
\end{proposition}

\begin{proof} In each SIN-group $X_n$ fix a base $\BB_n$ of open symmetric invariant neighborhoods of the neutral element $e$. Observe that any two neighborhoods $V_n\in\BB_n$, $V_m\in\BB_m$ with $n\le m$ commute:
\begin{equation}\label{COM}
V_nV_m=\bigcup_{x\in V_n}xV_m=\bigcup_{x\in V_n}V_mx=V_mV_n.
\end{equation}

Also for every $U\in\BB_n$  the $U^{\LR}$-ball centered at a point $x\in X$ has the form
$$B(x,U^{\LR})=\{y\in G_n:y\in xU\cap Ux\}=\{y\in G_n:y\in xU\}.$$
This implies that for any sequence of invariant neighborhoods $(U_n)_{n\in\w}\in\prod_{n\in\w}\BB_n$ we get
$$B(e;\sum_{n\ge 0}U_n^{\LR})=\LP\limits_{n\in\w}U_n\mbox{ where }\LP_{n\in\w}U_n=\bigcup_{n\in\w}\LP_{0\le n\le m}U_n\mbox{ and }\LP_{0\le n\le m}U_n=U_0U_1\cdots U_m.$$

Now, after this preparation, we are ready to prove that the uniform space $\ulim X_n$  is a topological group. First we prove that $\ulim X_n$ is a semitopological group, which means that the left and right shifts 
$$l_a:\ulim X_n\to\ulim X_n,\;\;l_a:x\mapsto ax\quad\mbox{and}\quad
r_a:\ulim X_n\to\ulim X_n,\;\;r_a:x\mapsto xa,$$
are continuous for any $a\in \ulim X_n$. Find $k\in\w$ such that $a\in X_k$  and observe that for every $m\ge k$ the restrictions $l_a|X_m:X_m\to\ulim X_n$ and $r_a|X_m:X_m\to\ulim X_n$ are uniformly continuous. Then the definition of $\ulim X_n$ implies that the shifts $l_a$ and $r_a$ are uniformly continuous and hence continuous. By the same reason, the inversion $$\ulim X_n\to\ulim X_n,\;\;x\mapsto x^{-1},$$is (uniformly) continuous.

It remains to prove that the group operation 
$$\ulim X_n\times\ulim X_n\to\ulim X_n,\quad (x,y)\mapsto xy,$$ is continuous at the neutral element $e$ of $X$.

Take any neighborhood $O_e\subset \ulim X_n$ of zero. By Theorem~\ref{topbase}, we can assume that $O_e$ is of the basic form 
$O_e=B(e;\sum_{n\ge 0}U^{\LR}_n)$ for some sequence  $(U_n)_{n\in\w}\in\prod_{n\in\w}\BB_n$.
The continuity of the binary operation in the SIN-groups $X_n$ yields a sequence of invariant neighborhoods $(V_n)_{n\in\w}\in\prod_{n\in\w}\BB_n$ such that $V_nV_n\subset U_n$ for every $n\in\w$. 

By induction on $m\in\w$ we shall prove that
\begin{equation}\label{EQ}
\big(\LP_{n\le m}V_n)\cdot\big(\LP_{n\le m}V_n\big)\subset \LP_{n\le m}U_n.
\end{equation}
For $m=0$ this follows from the choice of the neighborhood $V_0$. Assume that the inclusion (\ref{EQ}) has been proved for some $m=k$. For $m=k+1$ it also holds because:
$$
\begin{aligned}
\big(\LP_{n\le m}V_n)&\cdot\big(\LP_{n\le m}V_n\big)=
\big(\LP_{n<m}V_n\big)\cdot V_m\cdot \big(\LP_{n<m}V_n\big)\cdot V_m=\\
&\big(\LP_{n<m}V_n\big)\cdot\big(\LP_{n<m}V_n\big)\cdot V_mV_m\subset 
\big(\LP_{n<m}U_n\big)\cdot U_m=\LP_{n\le m}U_n.
\end{aligned}
$$
Here we have used the inductive assumption and the equality (\ref{COM}).

The sequence $(V_n)$ determines the neighborhood $B(e;\sum_{n\ge 0}V_n^{\LR})\subset\lim X_n$ witnessing the continuity of the group operation of $\ulim X_n$ at $e$:
$$
\begin{aligned}
B\big(e;\sum_{n\ge0}V_n^{\LR}\big)&\cdot B\big(e;\sum_{n\ge0}V_n^{\LR}\big)=\big(\LP_{n\in\w}V_n\big)\cdot
\big(\LP_{n\in\w}V_n\big)=\\
&\bigcup_{m\in\w}\big(\LP_{n\le m}V_n\big)\cdot
\big(\LP_{n\le m}V_n\big)\subset \bigcup_{m\in\w}\LP_{n\le m}U_n=\LP_{n\in\w}U_n=B\big(e,\sum_{n\ge 0}U_n^{\LR}\big).
\end{aligned}
$$
\end{proof}

\begin{remark} Proposition~\ref{group} identifies the direct limit $\glim X_n$ of a tower $(X_n)$ of SIN-groups in the category of topological groups with its direct limit $\ulim X_n$ in the category of uniform spaces. We do not know if such an identification still holds beyond the class of SIN-group. So, the problem of explicit description of the topology of the direct limits in th category of topological groups remains open. This problem was addressed in \cite{TSH}, \cite{HSTH}, \cite{Glo03}, \cite{Glo07}, \cite{BaR} where some partial answers are given.
\end{remark}

Now, let us switch to the direct limits in the category of locally convex spaces.
By the direct limit of a tower $(X_n)_{n\in\w}$ of (locally convex) linear topological spaces in the category of (locally convex) linear topological spaces we understand the union $X=\bigcup_{n\in\w}X_n$ endowed with the largest topology that turns $X$ into a (locally convex) linear topological space and making the identity operators $X_n\to X$ continuous.

\begin{proposition}\label{linear} 
The uniform direct limit $\ulim X_n$ of a tower $(X_n)_{n\in\w}$ of (locally convex) linear topological spaces is a (locally convex) linear topological space. Consequently, the identity map $\ulim X_n\to \llim X_n$ (as well as $\ulim X_n\to\lclim X_n$) is a homeomorphism.
\end{proposition}

\begin{proof} Observe that $X=\ulim X_n$, being the union of a tower of linear spaces,  is a linear space over a field $\IK$ (equal to the field of real or complex numbers). For each linear topological space $X_n$ consider the family $\mathcal O_n$ of open neighborhoods $U_n\subset X_n$ of zero such that $\lambda\cdot U\subset U$ for all $\lambda\in\IK$ with $|\lambda|\le1$. If the space $X_n$ is locally convex, we shall additionally assume that each set $U\in\mathcal O_n$ is convex.

By Proposition~\ref{group}, the uniform direct limit $\ulim X_n$ is a topological group with respect to the addition operation, and by Theorem~\ref{topbase} the family
$$\mathcal B_0=\Big\{\sum_{n\in\w}V_n:(V_n)_{n\in\w}\in\prod_{n\in\w}\mathcal O_n\Big\}$$is a neighborhood base of the topology on $\ulim X_n$ at zero.
It is clear that $\mathcal B_0$ has two properties:
\begin{itemize}
\item for any $U\in\mathcal B_0$ there is $V\in\mathcal B_0$ with $V+V\subset U$;
\item for any $U\in\mathcal B_0$ and $\lambda\in\IK$ with $|\lambda|\le 1$ we get $\lambda U\subset U$.
\end{itemize}

By \cite[1.2]{Sh}, those properties imply that $\ulim X_n$ is a topological linear space. 

If all the spaces $X_n$ are locally convex, then each set $U\in\mathcal B_0$, being the sum of convex sets, is convex, and hence the space $\ulim X_n$ is locally convex.
\end{proof}

As we already know, for certain tower $(X_n)_{n\in\w}$ of uniform spaces the identity map $\ulim X_n\to\tlim X_n$ is discontinuous. The following proposition detect towers $(X_n)$ for which that map is a homeomorphism. 

\begin{proposition} Let $(X_n)_{n\in\w}$ be a tower of uniform spaces. If each space $X_n$ is locally compact, then the identity map $\tlim X_n\to\ulim X_n$ is a homeomorphism.
\end{proposition}

\begin{proof} Since the map $\tlim X_n\to\ulim X_n$ always is continuous it suffices to check the continuity of the identity map $\ulim X_n\to\tlim X_n$ at each point $x\in\ulim X_n$. Pick any open neighborhood $O_t(x)$ of $x$ in $\tlim X_n$ and let $m=|x|$.  Since the space $X_{m}$ is locally compact, there is a closed entourage $U_{m}\in\U_{X_{m}}$ such that the ball $B(x,U_{m})$ is compact and lies in the open set $O(x)\cap X_{m}$ of $X_{|x|}$. Use the compactness of the set $B(x;U_m)\subset O(x)\cap X_{m+1}$ in the locally compact space $X_{m+1}$ in order to find a closed entourage $U_{m+1}\in\U_{X_{m+1}}$ such that $B(B(x;U_m);U_{m+1})=B(x;U_m+U_{m+1})$ is compact and lies in the open set $O_t(x)\cap X_{m+1}$.
Continuing by induction, construct a sequence $(U_n)_{n\ge n}\in\prod_{n\ge m}\U_{X_n}$ of closed entourages such that for every $n\in\w$ the ball $B(x;\sum_{m\le i\le n}U_i)$ is compact and lies in the open set $O_t(x)\cap X_n$. By Theorem~\ref{topbase}, the set $O_u(x)=B(x;\sum_{n\ge|x|}U_n)$ is a neighborhood of $x$ in $\ulim X_n$. Since $O_u(x)\subset O_t(x_0)$, the identity map $\ulim X_n\to\tlim X_n$ is continuous at the point $x$.
\end{proof}

Finally, we discuss the relation of the uniform direct limits to small box products.

By the small box-product of pointed topological spaces $(X_i,*_i)$, $i\in\I$, we understand the subspace
$$\cbox_{i\in\I}X_i=\big\{(x_i)_{i\in\I}\in\square_{i\in\I}X_i:\{i\in\I:x_i\ne *_i\}\mbox{ is finite}\big\}$$of the box-product $\square_{i\in\I}X_i$. The latter space in the Cartesian product $\prod_{i\in\I}X_i$ endowed with the box-topology generated by the boxes $\prod_{i\in\I}U_i$ where $U_i\subset X_i$, $i\in\I$, are open sets. 

If each space $X_i$, $i\in\I$, is uniform, then the box-product $\square_{i\in\I}X_i$ carries the box-uniformity generated by the entourages
$$\{\big((x_i),(y_i)\big)\in(\square_{i\in\I}X_i)^2:\forall i\in\I\;\;(x_i,y_i)\in U_i\big\}$$where $U_i\in\U_{X_i}$ for all $i\in\w$. The small box-product $\cbox_{i\in\I}X_i$ carries the uniformity inherited from $\square_{i\in\I}X_i$.

If the index set $\I$ is finite, then  the product $\cbox_{i\in\I}X_i=\square_{i\in\I}X_i$ turns into usual Tychonoff product $\prod_{i\in\I}X_i$ (endowed with the uniformity of Tychonoff product).

For any subset $\mathcal J\subset\I$ the small box-product $\cbox_{i\in\J} X_i$ will be identified with the subspace $$\{(x_i)_{i\in\I}\in\cbox_{i\in\I}X_i:\forall i\in\I\setminus\J\;(x_i=*_i)\}$$ of $\cbox_{i\in\I}X_i$.

So, for any sequence $X_n$, $n\in\w$, of pointed uniform spaces the small box-products $\cbox_{i\le n}X_i=\prod_{i\le n}X_i$, $n\in\w$, form a tower
$$X_0\subset X_0\times X_1\subset\cdots \cbox_{i\le n}X_i\subset\dots
$$ whose union coincides with $\cbox_{i\in\w}X_i$. 

Since the identity inclusions $\cbox_{i\le n}X_i\to\cbox_{i\in\w}X_i$, $n\in\w$, are uniformly continuous, the identity map $\id:\ulim \cbox_{i\le n}X_i\to\cbox_{i\in\w}X_i$ is uniformly continuous. 

\begin{proposition}\label{box} For any sequence $X_n$, $n\in\w$, of pointed uniform spaces, the identity map $\id:\cbox_{n\in\w}X_n\to\ulim \cbox_{i\le n}X_i$ is a homeomorphism.
\end{proposition}

\begin{proof} It suffices to check the continuity of the identity map $\id:\cbox_{n\in\w}X_n\to\ulim \cbox_{i\le n}X_i$ at arbitrary point $\vec x\in\cbox_{i\in\w}X_i$. Let $k=\min\{n\in\w:\vec x\in\cbox_{i\le n}X_i\}$.

Given a neighborhood $O_u(\vec x)\subset \ulim\cbox_{i\le n}X_i$ of $\vec x$ we need to construct a neighborhood $O_b(\vec x)\subset\cbox_{i\in\w}X_i$ such that $O_b(\vec x)\subset O_u(\vec x)$.  By Theorem~\ref{topbase}, we can assume that $O_u(\vec x)$ is of the form
$O_u(\vec x)=B(\vec x,\sum_{n\ge k}U_n)$ for some sequence 
$(U_n)_{n\ge k}$ of entourages $U_n$ in the spaces $\cbox_{i\le n}X_i=\prod_{i\le n}X_i$, $n\ge k$. We can also assume that each entourage $U_n$ is of the basic form $$U_n=\prod_{i\le n}U_{i,n}=\{\big((x_i),(y_i)\big)\in(\cbox_{i\le n}X_i)^2:\forall i\le n\;\;(x_i,y_i)\in U_{i,n}\}$$ for some open entourages $U_{i,n}\in \U_{X_i}$, $i\le n$.

It follows that $V_{k}=B(\vec x,U_k)$ is an open neighborhood of $\vec x$ in $\cbox_{i\le k}X_i$ and $O_b(\vec x)=V_k\times\cbox_{i>k}B(*_i,U_{i,i})$ is an open neighborhood of $\vec x$ in $\cbox_{i\in\w}X_i$ such that $O_b(\vec x)\subset O_u(\vec x)$.
\end{proof}

\end{document}